\documentclass[times,sort&compress,3p]{elsarticle}
\journal{Journal of Multivariate Analysis}
\usepackage[labelfont=bf]{caption}

\usepackage{amsmath,amsfonts,amssymb,amsthm,booktabs,color,epsfig,graphicx,hyperref,url}

\newtheorem{theorem}{Theorem}

\newtheorem{corollary}{Corollary}
\newtheorem{lemma}{Lemma}

\usepackage{caption}
\usepackage{here}
\usepackage{subcaption}
\renewenvironment{proof}{{\bfseries Proof.}}{\qed}
\begin{document}

\begin{frontmatter}

\title{Heterogeneous hypergeometric functions with two matrix arguments 
and the exact distribution of the largest eigenvalue of a singular beta-Wishart matrix}

\author[A1]{Koki Shimizu}
\author[A1]{Hiroki Hashiguchi}

\address[A1]{Tokyo University of Science, 1-3 Kagurazaka, Shinjuku-ku, Tokyo 162-8601,Japan}

\cortext[mycorrespondingauthor]{Corresponding author. Email address: \url{1420702@ed.tus.ac.jp}~(K. Shimizu).}

\begin{abstract}
This paper discusses certain properties of heterogeneous hypergeometric functions with two matrix arguments.
These functions are newly defined but have already appeared in statistical literature and are useful when dealing with the derivation of certain distributions 
for the eigenvalues of singular beta-Wishart matrices. The joint density function of the eigenvalues and the distribution 
of the largest eigenvalue can be expressed in terms of heterogeneous hypergeometric functions.
Exact computation of the distribution of the largest eigenvalue is conducted for a real and complex cases. 
\end{abstract}


\end{frontmatter}
\section{Introduction\label{sec:1}}
Hypergeometric functions with one or two matrix arguments appear in distribution theory in multivariate analysis. 
The density and distribution functions of eigenvalues of a non-singular real Wishart matrix were given by
James~\cite{James1964} and Constantine~\cite{Constantaine1963} in terms of hypergeometric functions. 
 The exact distributions of the largest and smallest eigenvalues were derived by Sugiyama~\cite{Sugiyama1967} and Khatri~\cite{Khatri1972}, respectively. 
Ratnarajah et al.~\cite{Ratnarajah2005c} extended these results to the complex case and applied them to the channel capacity of wireless communication systems. 
Recently, these classes have been generalized to a beta-Wishart ensemble or non-singular beta-Wishart matrix that includes a classical matrix of real, complex, and quaternion cases.
Koev and Edelman~\cite{Edelman2014} derived exact distributions of the extreme eigenvalues for a non-singular beta-Wishart matrix. 
The essence of this work lies in the use of Jack polynomials instead of zonal polynomials. 
D\'iaz-Garc\'ia and Guti\'errez-J\'aimez~\cite{Garcia2011} discussed Wishart matrices on a real finite-dimensional normed division algebra that coincides with real numbers, complex numbers, quaternions, and octonions. 
In situations of normed division algebra, the Stiefel manifold and Haar measure can be defined as well as in the classical real case.
D\'iaz-Garc\'ia~\cite{Garcia2013} derived a useful integral formula for Jack polynomials over the Steifel manifold on a normed division algebra. 

The singular real Wishart matrix was defined by Mitra~\cite{Mitra1970}, Khatri~\cite{Khatri1970}, Srivastava and Khatri~\cite{Srivastava1979} among others.
Uhlig~\cite{{Uhlig1994}} gave the density of a singular real Wishart matrix and discussed the need for such distributions in his Bayesian analysis. 
The properties of this distribution were well studied by D\'iaz-Garc\'ia et al.~\cite{Garcia1997}, Srivastava~\cite{Srivastava2003}, Bodnar and Okhrin~\cite{Bodnar2008}, Bodnar et al.~\cite{Bodnar2018, Bodnar2014}.
Ratnarajah and Vaillancourt~\cite{Ratnarajah2005a,Ratnarajah2005b} derived some results for a singular complex Wishart matrix, and
Li and Xue~\cite{Li2010} studied a singular quaternion Wishart matrix. 
However, the exact distributions of the extreme eigenvalues have not been derived in singular real, complex, and quaternion cases. 

In this paper, we propose the key ideas for developing the distribution theory of eigenvalues for a singular beta-Wishart matrix. 
We show that the eigenvalue distributions of a singular beta-Wishart matrix are expressed 
in terms of heterogeneous hypergeometric functions with two matrix arguments. 
In the derivation of the heterogeneous hypergeometric functions, zero eigenvalues are substituted into 
the parts of arguments of zonal or Jack symmetric polynomials. 
In Section~\ref{sec:02}, we provide some definitions and preliminaries that include 
classical hypergeometric functions in terms of zonal polynomials. 
We define the heterogeneous hypergeometric functions with two matrix arguments 
and derive some useful formulas for them. 
Furthermore, we show that the integral over the Stiefel manifold can be converted into an integral over the orthogonal group.
In Section~\ref{sec:03}, we derive the exact distribution of the largest eigenvalue of a singular real Wishart matrix.  
The distribution of the largest eigenvalue is applied using the formula of Sugiyama~\cite{Sugiyama1967}. 
We then conduct numerical experiments for a theoretical distribution. 
Numerical results are calculated using the algorithm of Hashiguchi et al.~\cite{Hashiguchi2000} for zonal polynomials. 
In Section~\ref{sec:04}, we discuss the real finite division algebra and define the heterogeneous hypergeometric functions of parameter $\beta>0$. 
We restrict the parameter $\beta=1$~(real case), $2$~(complex case), $4$~(quaternion case) and present some useful formulas for these values on a real finite division algebra. 
In Section~\ref{sec:05}, we derive distributions of a singular beta-Wishart matrix.
 We also illustrate numerical computation for $\beta =2$ and application to the channel capacity of a multiple input and single output system (MISO).
Using the heterogeneous hypergeometric functions, Shimizu and Hashiguchi~\cite{Shimizu2020} discussed a singular $F$-matrix that is important to inference in statistical multivariate analysis for $\beta=1$, for example, in MANOVA, classification and test for equality of covariance matrices. 
\section{Heterogeneous hypergeometric functions ${_pF_q^{(m,n)}}$}
\label{sec:02}
In this section, we consider only the case of real numbers and define heterogeneous hypergeometric functions with two matrix arguments. These often appear in the density functions of random matrices for a singular Wishart matrix. Ratnarajah and Vaillancourt~\cite{Ratnarajah2005a,Ratnarajah2005b} used such hypergeometric functions in the derivation of the density function  of a singular complex Wishart matrix. 
The exterior product for real matrix $X$ is written by $(dX)$ as defined in Muirhead~\cite{Muirhead1982} and Gupta and Nagar~\cite{Gupta1999}. 
The set of all $m\times n$ matrices $H_1$ with orthonormal columns is called the Stiefel manifold, denoted by $V_{n,m}$ where $n\leq m$, 
\begin{align}
\label{stiefel:def}
V_{n,m}=\{ H_1\mid H_1^\top H_1=I_n\}. %スティーフェルの定義
\end{align}
We note that $V_{m,m}=O(m)$, where $O(m)$ is the orthogonal group of order $m$. 
If $H_1\in V_{n,m}$, then we define 
\begin{align*}
(H_1^\top dH_1)=\bigwedge_{i=1}^{n}\bigwedge_{j=i+1}^{m}h_j^\top dh_i,
\end{align*}
where $H=(H_1:H_2)\in O(m)$. 
The volumes of $V_{n,m}$ and $O(m)$ are given by
\begin{align*}
\mathrm{Vol}(V_{n,m})=\int_{H_1\in V_{n,m}}(H_1^\top dH_1)=\frac{2^n\pi^{mn/2}}{\Gamma_n(m/2)},\quad 
\mathrm{Vol}(O(m))=\int_{H\in O(m)}(H^\top dH)=\frac{2^m\pi^{m^2/2}}{\Gamma_m(m/2)},
\end{align*}
respectively, where the multivariate gamma function is
\begin{align*}
\Gamma_m(c)=\pi^{m(m-1)/4}\prod_{i=1}^{m}\Gamma \biggl(c-\frac{i-1}{2}\biggl), \ \ \mathrm{Re}(c)>n-1. %ガンマ関数の定義
\end{align*}
The differential form $(dH_1)$, defined as 
\begin{align*}
(dH_1)=\frac{(H_1^\top dH_1)}{\mathrm{Vol}(V_{n,m})}=\frac{\Gamma_n(m/2)}{2^n\pi^{mn/2}}(H_1^\top dH_1)
\end{align*}
is normalized as 
\begin{align*}
\int_{H_1\in V_{n,m}}(dH_1)=1.
\end{align*}
This differential form represents the Haar invariant probability measure on $V_{n,m}$. 

 For a positive integer $k$, let $\kappa=(\kappa_1,\kappa_2,\dots,\kappa_m)$ denote a partition of $k$ with $\kappa_1\geq\cdots\geq \kappa_m\geq 0$ and $\kappa_1+\cdots +\kappa_m=k$. The set of all partitions with lengths not longer than $m$ is denoted by $P^k_{m}=\{\ \kappa=(\kappa_1,\dots,\kappa_m)\mid \kappa_1+\dots+\kappa_m=k, \kappa_1\geq \kappa_2\geq\cdots\geq  \kappa_m \geq 0 \}$. 
The  Pochhammer symbol for a partition $\kappa$ is defined as $(\alpha)_\kappa=\prod_{i=1}^{m}\{\alpha-(i-1)/2\}_{\kappa_i}$, where $(\alpha)_k=\alpha(\alpha+1)\cdots (\alpha+k-1)$ and $(\alpha)_0=1$. 
Furthermore, let $X$ be an $m\times m$ symmetric matrix with eigenvalues $x_1,\dots,x_m$. 
Then, the zonal polynomial $C_\kappa(X)$ is defined as a symmetric polynomial in $x_1,\dots,x_m$. 
Chapter~7 of Muirhead~\cite{Muirhead1982} presents a detail definition. If the length of $\kappa$ is $m$ and the $m\times m$ symmetric matrix $X$ has a rank of $n\leq m$ so that $x_{n+1}=\cdots=x_m=0$, then we have 
\begin{align}
\label{eq:eq01}
C_{\kappa}(X)=
\begin{cases}
\quad 0, \quad \quad \quad \quad \quad \text{\ $n<m$,}\\
\quad C_{\kappa}(X_1), \quad \quad \quad \text{$n\geq m$},
\end{cases}
\end{align}
where $X_1=\mathrm{diag}(x_1,\dots,x_n)$. 
For integers $p,q\geq 0$ and $m\times m$ real symmetric matrices $A$ and $B$, the hypergeometric function with two matrices is defined as
\begin{align}
\label{eq:eq02}
{}_pF_q{}^{(m)}(\mbox{\boldmath$\alpha$};\mbox{\boldmath$\beta$};A,B) = \sum_{k=0}^\infty \sum_{\kappa \in P^k_m}\frac{(\alpha_1)_{\kappa} \cdots (\alpha_p)}{(\beta_1)_{\kappa} \cdots (\beta_q)}\frac{C_\kappa(A)C_\kappa(B)}{k!C_\kappa(I_m)}, 
\end{align}
where $\mbox{\boldmath$\alpha$}=(\alpha_1,\ldots,\alpha_p)^\top$, $\mbox{\boldmath$\beta$}=(\beta_1,\ldots,\beta_q)^\top$. For one matrix argument $A$, we also define ${}_pF_q{}(\mbox{\boldmath$\alpha$};\mbox{\boldmath$\beta$};A)$ as
\begin{align}
\label{eq:eq03}
{}_pF_q{}(\mbox{\boldmath$\alpha$};\mbox{\boldmath$\beta$};A)={}_pF_q{}^{(m)}(\mbox{\boldmath$\alpha$};\mbox{\boldmath$\beta$};A,I_m).
\end{align}
Then, the following relationship between (\ref{eq:eq02}) and (\ref{eq:eq03}) holds, 
\begin{align}
\label{eq:eq04}
{}_pF_q
{}^{(m)}(\mbox{\boldmath$\alpha$};\mbox{\boldmath$\beta$};A,B)=\int_{H\in O(m)}{_pF_q}(\mbox{\boldmath$\alpha$};\mbox{\boldmath$\beta$};AHBH^\top)(dH). 
\end{align}
%%%%%%%%%%%%%%%%%%%% %Heterogeneous hypergeometric Functions
To discuss the density function for a singular Wishart matrix, we define the heterogeneous hypergeometric functions as follows.
For an $m\times m$ symmetric matrix $A$ and an $n\times n$ symmetric matrix $B$, the heterogeneous hypergeometric functions are defined as 
\begin{align}
\label{eq:eq05}
{_pF_q}^{(s,t)}({\mbox{\boldmath$\alpha$}};\mbox{\boldmath$\beta$};A,B)=\sum_{k=0}^{\infty}\sum_{\kappa \in P^k_t}\frac{(\alpha_1)_\kappa \cdots(\alpha_p)_\kappa}{(\beta_1)_\kappa \cdots(\beta_q)_\kappa}\frac{C_\kappa(A)C_\kappa(B)}{k!C_\kappa(I_s)}, 
\end{align}
where $s=\mathrm{max}(m,n)$, $t=\mathrm{min}(m,n)$. 

We consider the case of $m\geq n$ as well as $s=m$ and $t=n$ hereafter. 
So we write ${_pF_q}^{(m,n)}$ instead of (\ref{eq:eq05}). If $m>n$, it is clear that
\begin{align*}
{_pF_q}^{(m,n)}({\mbox{\boldmath$\alpha$}};\mbox{\boldmath$\beta$};I_m,B)={_pF_q}({\mbox{\boldmath$\alpha$}};\mbox{\boldmath$\beta$};B) 
\end{align*}
and
\begin{align*}
{_pF_q}^{(m,n)}({\mbox{\boldmath$\alpha$}};\mbox{\boldmath$\beta$};A,I_n) &\neq{_pF_q}({\mbox{\boldmath$\alpha$}};\mbox{\boldmath$\beta$};A). 
\end{align*}
For an $m \times m$ matrix $B_1=\left( 
\begin{array}{cc}
B & O\\
O & O
\end{array}
\right)$, we have the following relationship from (\ref{eq:eq01}):
\begin{align*}
{_pF_q}^{(m,n)}({\mbox{\boldmath$\alpha$}};\mbox{\boldmath$\beta$};A,B)={_pF_q}^{(m)}({\mbox{\boldmath$\alpha$}};\mbox{\boldmath$\beta$};A,B_1).
\end{align*}

%%%%%%%%%%%%%%%%%%%%%%%%%%%%%%%%%%%%%
\begin{lemma}
For an $m\times m$ positive symmetric matrix $A$ and $m\times m$ symmetric positive semi-definite matrix $B$, we have
\begin{align*}
\int_{H\in O(m)}C_\kappa(AHBH^\top)(dH)=\frac{C_\kappa(A)C_\kappa(B)}{C_\kappa(I_m)}.
\end{align*}
\end{lemma}
\noindent
\begin{proof}
From the fundamental properties of zonal polynomials $C_\kappa(AB)=C_\kappa(A^{1/2}BA^{1/2})=C_\kappa(BA)$ and $C_\kappa(H^\top AH)=C_\kappa(A)$ for any $H\in O(m)$, 
\begin{align*}
\int_{H\in O(m)}C_\kappa(AHBH^\top)(dH)
&=\int_{H\in O(m)}C_\kappa(BH^\top AH)(dH). 
\end{align*}
Let $f_\kappa(A)$ be 
\begin{align*}
f_\kappa(A)=\int_{H\in O(m)}C_\kappa(BH^\top AH)(dH).
\end{align*}
From the proof of Theorem 7.2.5 in Muirhead~\cite{Muirhead1982}, $f_\kappa(A)$ must be a multiple of the zonal polynomial $C_\kappa(A)$; that is, $f_\kappa(A)=\lambda_\kappa C_\kappa(A)$.
Putting $A=I_m$ and using $f_\kappa(I_m)= C_\kappa(B)$, we have $\lambda_\kappa=\frac{C_\kappa(B)}{C_\kappa(I_m)}$. 
\end{proof}
%%%%%%%%%%%%%%%%%%%%%%%%%%%%%%%%%%%%%
\begin{theorem}
For an $m\times m$ positive symmetric matrix $A$ and an $n\times n$ symmetric matrix $B$, we have
\begin{align}
\label{splitting}
\int_{H_1\in V_{m,n}}C_\kappa(AH_1BH_1^\top)(dH_1)=\frac{C_\kappa(A)C_\kappa(B)}{C_\kappa(I_m)}.
\end{align}
\end{theorem}
\noindent
\begin{proof}
We refer to the proof of Lemma 9.5.3 of Muirhead~\cite{Muirhead1982}. 
For any $m\times(m-n)$ matrix $G$ with orthonormal columns that are orthogonal to those of $H_1$, $K\in O(m-n)$, $H=(H_1,H_2)\in O(m)$, and $m \times m$ matrix $B_1=\left( 
\begin{array}{cc}
B & O\\
O & O
\end{array}
\right)$, and we have $H_2=GK$ and 
\begin{align*}
AHB_1H^\top
=A
(H_1, GK)
\left( 
\begin{array}{cc}
B & O\\
O & O
\end{array}
\right) 
\left(
\begin{array}{c}
H_1^\top \\
KG^\top
\end{array}
\right) 
=A 
(H_1 B , O)
\left(
\begin{array}{c}
H_1^\top \\
K G^\top
\end{array}
\right) 
=AH_1 B H_1^\top.
\end{align*}
Because of $dh_{n+j}=Gdk_j, i\in\{1,\ldots, m-n\}$, it is clear that
\begin{align}
\label{eq:eq31}
(H^\top dH)&=(H_1^\top dH_1)(K^\top dK).
\end{align}
The volume $\mathrm{Vol}(V_{n,m})$ is written by
\begin{align}
\label{eq:eq07}
\mathrm{Vol}(V_{n,m})=\frac{\mathrm{Vol}(O(m))}{\mathrm{Vol}(O(m-n))}.
\end{align}
From (\ref{eq:eq31}) and (\ref{eq:eq07}), we have $(dH)=(dH_1)(dK)$, which means that $\displaystyle \int_{O(m-n)}(dK)=1$. 

Therefore, 
\begin{align*}
\int_{H_1\in V_{n,m}}C_\kappa(AH_1BH_1^\top)(dH_1)
&=\int_{K\in O(m-n)}(dK)\int_{H_1\in V_{n,m}}C_\kappa(AH_1BH_1^\top)(dH_1)\\
&=\int_{K\in O(m-n)}\int_{H_1\in V_{n,m}}C_\kappa(AH_1BH_1^\top)(dH_1)(dK)\\
&=\int_{H\in O(m)}C_\kappa(AHB_1H^\top)(dH)
=\frac{C_\kappa(A)C_\kappa(B)}{C_\kappa(I_m)}.
\end{align*}
\end{proof}
%%%%%%%%%%%%%%%%%%%%

Using (\ref{eq:eq01}), the following relationship holds from Lemma~1 and Theorem~1. 
\begin{align*}
\int_{H_1\in V_{n,m}}C_\kappa(AH_1BH_1^\top)(dH_1)
&=\int_{H\in O(m)}C_\kappa(AHB_1H^\top)(dH).
\end{align*}
The next equation (\ref{eq:eq35}) is immediately  obtained from the above result. 
\begin{align}
\label{eq:eq35}
\displaystyle \int_{H_1\in V_{n,m}}{_pF_q}(\mbox{\boldmath$\alpha$};\mbox{\boldmath$\beta$};AH_1BH_1^\top)(dH_1)&=\displaystyle \int_{H\in O(m)}{_pF_q}(\mbox{\boldmath$\alpha$};\mbox{\boldmath$\beta$};AHB_1H^\top)(dH).
 \end{align}

%%%%%%%%%%%%%%%%%%%%%%%%%%%%%%%%%%%
The following Corollary~1 provides a result similar to that in (\ref{eq:eq04}). 
\begin{corollary}
For an $m\times m$ symmetric matrix $A$, $B=\mathrm{diag}(b_1,\dots,b_n)$ and $m\times m$ diagonal matrix $B_1=\mathrm{diag}(b_1,\dots,b_n,0,\dots,0)$ where $m\geq n$, we have
\begin{align*}
{_pF_q}^{(m,n)}(\mbox{\boldmath$\alpha$};\mbox{\boldmath$\beta$};A,B)
&=\displaystyle \int_{H\in O(m)}{_pF_q}(\mbox{\boldmath$\alpha$};\mbox{\boldmath$\beta$};AHB_1H^\top)(dH)
=\displaystyle \int_{H_1\in V_{n,m}}{_pF_q}(\mbox{\boldmath$\alpha$};\mbox{\boldmath$\beta$};AH_1BH_1^\top)(dH_1).
 \end{align*}
 \end{corollary}
 \noindent
\begin{proof}
From (\ref{eq:eq35}), we have 
\begin{align*}
{_pF_q}^{(m,n)}(\mbox{\boldmath$\alpha$};\mbox{\boldmath$\beta$};A,B)&={_pF_q}^{(m)}(\mbox{\boldmath$\alpha$};\mbox{\boldmath$\beta$};A,B_1)
=\displaystyle \int_{H\in O(m)}{_pF_q}(\mbox{\boldmath$\alpha$};\mbox{\boldmath$\beta$};AHB_1H^\top)(dH)
=\displaystyle \int_{H_1\in V_{n,m}}{_pF_q}(\mbox{\boldmath$\alpha$};\mbox{\boldmath$\beta$};AH_1BH_1^\top)(dH_1).
\end{align*}
\end{proof}
%%%%%%%%%%%%

Two particular cases of (\ref{eq:eq05}) are listed as 
\begin{align}
\label{eq:eq08}
{_0F_0}^{(m,n)}(A,B)&=\sum_{k=0}^{\infty}\sum_{\kappa \in P^k_n}\frac{C_\kappa(A)C_\kappa(B)}{k!C_\kappa(I_m)},
\quad
{_1F_0}^{(m,n)}(\alpha_1;A,B)=\sum_{k=0}^{\infty}\sum_{\kappa \in P^k_n}\frac{(\alpha_1)_\kappa}{k!}\frac{C_\kappa(A)C_\kappa(B)}{C_\kappa(I_m)}
\end{align}
and from Corollary~1, the related integral formulas of (\ref{eq:eq08}) are given in Corollary~2. 
%%%%%%%%%%%%%%%%%%%%%%%%%%%%%%%%%%%%%%%%
\begin{corollary}
For an $m\times m$ symmetric matrix $A$,  $B=\mathrm{diag}(b_1,\dots,b_n)$ and a nonnegative integer $r$, we have
\begin{align}
\label{0F0integral}
 {}_0F_0{}^{(m,n)}\left(A,B\right)&=\displaystyle \int_{H_1\in V_{n,m}}\mathrm{etr} \left(AH_1BH_1^\top \right)(dH_1),\\
 \label{eq:eq11}
 {_1F_0}^{(m,n)}(\alpha_1;A,B)&=\int_{H_1\in V_{n,m}}|I_m-AH_1BH_1^\top|^{-\alpha_1}(dH_1),
 \end{align}
 where $\mathrm{etr}(\cdot)=\mathrm{exp}(\mathrm{tr}(\cdot))$ and $|A|$ is the determinant of the matrix $A$.
  \end{corollary}
  \noindent
 \begin{proof}
For an $m \times m$ diagonal matrix, $B_1=\mathrm{diag}(b_1,\dots,b_n,0,\dots,0)$ and $(H_1,H_2)\in H$, where $H_1\in V_{n,m}$. The special cases of (\ref{eq:eq03}) are represented as ${}_0F_0{}\left(A\right)=\mathrm{etr}\left(A\right)$ and ${}_1F_0{}\left(\alpha_1;A\right)=|I_m-A|^{-\alpha_1}$. Then
\begin{align*}
{}_0F_0{}^{(m,n)}\left(A,B\right)
&={}_0F_0{}^{(m)}\left(A,B_1\right)
=\int_{H\in O(m)}{}_0F_0{}\left(AHB_1H^\top\right)(dH)\\
&=\displaystyle \int_{H\in O(m)}\mathrm{etr} \left(AHB_1H^\top\right)(dH)
=\displaystyle \int_{H_1\in V_{n,m}}\mathrm{etr} \left(AH_1BH_1^\top\right)(dH_1).
\end{align*}
The identity (\ref{eq:eq11}) can also be derived in the same way as above. 
 \end{proof}
%%%%%%%%%%%%%%%%%%%%%%%%%%%%%%%%%%%%%%%%%%
\begin{corollary}
For an $m\times m$ symmetric matrix $A$ and $B=\mathrm{diag}(b_1,\dots,b_n)$, we have
\begin{align*}
{}_0F_0{}^{(m,n)}\left(A+I_m,B\right)=\mathrm{etr}(B)~{}_0F_0{}^{(m,n)}\left(A,B\right).
\end{align*}
\end{corollary}
\noindent
\begin{proof}
For an $m \times m$ diagonal matrix $B_1=\mathrm{diag}(b_1,\dots,b_n,0,\dots,0)$ and $(H_1,H_2)\in H$ where $H_1\in V_{n,m}$, 
\begin{align*}
{}_0F_0{}^{(m,n)}\left(I_m+A,B\right)&=\int_{H_1\in V_{n,m}}{}_0F_0{}^{(m,n)}\left\{(I_m+A)H_1BH_1^\top\right\}(dH_1)
=\int_{H\in O(m)}{}_0F_0{}\left\{(I_m+A)HB_1H^\top\right\}(dH)\\
&=\mathrm{etr}(B)\int_{H\in O(m)}\mathrm{etr}\left(AHB_1H^\top\right)(dH)=\mathrm{etr}(B)~{}_0F_0{}^{(m)}(A,B_1)=\mathrm{etr}(B)~{}_0F_0{}^{(m,n)}(A,B).
\end{align*}
\end{proof}
\section{Exact distribution of the largest eigenvalue of a singular Wishart matrix}
\label{sec:03}
Suppose that an $m\times n$ real Gaussian random matrix $X$ is distributed as $X\sim$ $\mathcal{N}_{m,n}(O,\Sigma \otimes I_n)$, 
where $O$ is the $m \times n$ zero matrix, $\Sigma > 0$, and $\otimes$ is the Kronecker product. This means that 
the column vectors of $X$ are an \rm{i.i.d.} sample of size $n$ from $\mathcal{N}_{m}(0, \Sigma)$, where $0$ is the $m$-dimensional zero vector.
%$m$ dimensional normal distribution with zero mean vector and covariance matrix $\Sigma$. 
The non-zero eigenvalues of $\Sigma$ are denoted by $\lambda_1, \lambda_2, \dots, \lambda_m$, 
where $\lambda_1\geq \lambda_2 \geq \cdots \geq \lambda_m>0$. 
Then the random matrix $W=XX^\top $ is called a non-singular real Wishart matrix. 
The eigenvalues of $W$ are denoted by $\ell_1,\dots,\ell_m$, with $\ell_1>\ell_2>\cdots>\ell_m>0$. 
If $n<m$, then $W$ is said to be a singular real Wishart matrix. 
The first $n$ eigenvalues are not zero and the remaining $m-n$ eigenvalues $\ell_{n+1},\dots,\ell_m$ are all zero. 
Uhlig~\cite{Uhlig1994} derived that the density function of $W$ for a singular case is given as
\begin{align*}
f(W)=\frac{\pi^{(-mn+n^2)/2}}{2^{mn/2}\Gamma_n(n/2)(\mathrm{det}\Sigma)^{n/2}}\mathrm{etr}(-\Sigma^{-1} W/2)(\mathrm{det}L_1)^{(n-m-1)/2},
\end{align*}
where $L_1=\mathrm{diag}(\ell_1,\dots,\ell_n)$. 
Srivastava~\cite{Srivastava2003} showed that the joint density function of $\ell_1, \dots, \ell_n$ is given as
\begin{align}
\label{eq:eq12}
 C\left(\prod_{i=1}^n \ell_i^{(m-n-1)/2}\right)\left(\prod_{i<j}^{n}(\ell_i - \ell_j) \right) \int_{H_1\in V_{n,m}}\mathrm{etr} \left(-\frac{1}{2}\Sigma^{-1}H_1L_1H_1^\top \right), 
\end{align}
 where $C=\frac{2^{-nm/2}\pi^{n^2/2}|\Sigma|^{-n/2}}{\Gamma_n(\frac{n}{2})\Gamma_n(\frac{m}{2})}$. From (\ref{0F0integral}), the equation (\ref{eq:eq12}) is also expressed as 
 \begin{align}
\label{eq:eq13}
 C\left(\prod_{i=1}^n \ell_i^{(m-n-1)/2}\right)\left(\prod_{i<j}^{n}(\ell_i - \ell_j) \right)
 {}_0F_0{}^{(m,n)}\left(-\frac{1}{2}\Sigma^{-1},L_1\right).
\end{align}
%%%%%%%%%%%%%%%%%%%%%%%%%%%%%%%%%%%%%%%%%%%%%%%%%%%%
The following two lemma's are required in order to integrate (\ref{eq:eq13}) with respect to $\ell_2, \dots, \ell_n$.
\begin{lemma} \label{lem-Sugiyama01}
Let $L=\mathrm{diag}(\ell_1,\dots ,\ell_n)$ and the length of $\kappa$ be equal to or less than n. 
Then the following equation holds: 
\begin{align*}
&&\int_{1>\ell_1>\ell_2>\cdots \ell_n>0}|L|^{t-(n+1)/2}|I_n-L|^{u-(n+1)/2}C_\kappa(L)\prod_{i<j}^{n}(\ell_i-\ell_j)\prod_{i=1}^{n}d\ell_i =\frac{\Gamma_n(n/2)\Gamma_n(t,\kappa)\Gamma_n(u)C_\kappa(I_n)}{\pi^{n^2/2}\Gamma_n(t+u,\kappa)}, 
\end{align*}
where $\mathrm{Re}(t)>(n-1)/2, \mathrm{Re}(u)>(n-1/2), \Gamma_n(\alpha,\kappa)=(\alpha)_\kappa\Gamma(\alpha)$. 
\end{lemma}
%%%%%%%%%%%%%%%%%%%%%%%%%%%%%%%%%%%%%%%%%%%%%%%
\begin{lemma} \label{lem-Sugiyama02}
Let $X_1=\mathrm{diag}(1,x_2,\dots ,x_n)$ and $X_2=\mathrm{diag}(x_2,\dots,x_n)$ with $x_2>\cdots >x_n>0$. Then the following equation holds: 
\begin{align*}
&&\int_{1>x_2>\cdots x_n>0}|X_2|^{t-(n+1)/2}C_\kappa(X_1)\prod_{i=2}^{n}(1-x_i)\prod_{i<j}(x_i-x_j)\prod_{i=1}^{n}dx_i=(nt+k)(\Gamma_n(n/2)/\pi^{n^2/2})\frac{\Gamma_n(t,\kappa)\Gamma_n\{(n+1)/2\}C_\kappa(I_n)}{\Gamma_n\{t+(n+1)/2,\kappa\}}. 
\end{align*}
\end{lemma}
%%%%%%%%%%%%%%%%%%%%%%%%%%%%%%%%%%%%%%%%%%%%%%

From Lemma~\ref{lem-Sugiyama01}, Sugiyama~\cite{Sugiyama1967} derived Lemma~\ref{lem-Sugiyama02} 
for the derivation of the exact distribution of $\ell_1$ for a non-singular real Wishart matrix. 
Shinozaki et al.~\cite{Shinozaki2018} gave the distribution of the largest eigenvalue under an elliptical population 
using Lemmas~\ref{lem-Sugiyama01} and \ref{lem-Sugiyama02}. 
In the case of a singular Wishart matrix, the exact distribution of $\ell_1$ of $W$ is given in Theorem~\ref{th:dis-l1}. 

\begin{theorem} \label{th:dis-l1}
Let $W\sim \mathcal{W}_m(n,\Sigma)$, where $m>n$. 
Then the distribution function of the largest eigenvalue $\ell_1$ of $W$ is given as
\begin{align*}
\Pr (\ell_1<x)
=\frac{\Gamma_n\{(n+1)/2\}(\frac{x}{2})^{nm/2}}{\Gamma_n\{(n+m+1)/2\}|\Sigma|^{n/2}}
\; {{}_1F_1{}}^{(m,n)}\left(\frac{m}{2};\frac{n+m+1}{2};-\frac{1}{2}x\Sigma^{-1},I_n\right ). 
\end{align*}
\end{theorem}
\noindent
\begin{proof}
The joint density of $\ell_1, \ell_2, \dots, \ell_n$ is given in (\ref{eq:eq13}) as
\begin{align*}
f(\ell_1,\dots,\ell_n)&= C(\mathrm{det}L_1)^{(m-n-1)/2}\prod_{i<j}^{n}(\ell_i - \ell_j){}_0F_0{}^{(m,n)}\left(-\frac{1}{2}\Sigma^{-1},L_1\right)\\
&=C(\mathrm{det}L_1)^{(m-n-1)/2}\prod_{i<j}^{n}(\ell_i-\ell_j)\sum_{k=0}^{\infty}\sum_{\kappa \in P^k_n}\frac{C_\kappa(-\frac{1}{2}\Sigma^{-1})C_\kappa(L_1)}{k!C_\kappa( I_m)}.
\end{align*}
Translating $\ell_i$ to $x_i=\ell_i/\ell_1$,  $i\in\{2,\ldots, n\}$, and using Lemma~\ref{lem-Sugiyama02}, the density function of $\ell_1$ is given as
\begin{align*}
f(\ell_1)&=\mathrm{C}\sum_{k=0}^{\infty}\sum_{\kappa \in P^k_n}\int_{1>x_2>\cdots x_n>0}(\mathrm{det}X_2)^{(m-n-1)/2}\prod_{i=2}^{n}(1-x_i)\prod_{2\leq i<j}(x_i-x_j)C_\kappa(X_1)\ell_1^{mn/2+k-1} \frac{C_\kappa(-\frac{1}{2}\Sigma^{-1})}{k!C_\kappa( I_m)}\\
&=C (nm/2+k) \Gamma_n(n/2)\sum_{k=0}^{\infty}\sum_{\kappa \in P^k_n} \frac{\Gamma_n(m/2,\kappa)\Gamma_n\{(n+1)/2\}C_\kappa(I_n)}{\pi^{n^2/2}\Gamma_n\{n+m+1)/2,\kappa \}} \ell_1^{mn/2+k-1}\frac{C_\kappa(-\frac{1}{2}\Sigma^{-1})}{k!C_\kappa( I_m)},
\end{align*}
where $X_1=\mathrm{diag}(1,x_2,\dots,x_n)$ and $X_2=\mathrm{diag}(x_2,\dots,x_n), x_2>\cdots >x_n>0$. 
Moreover, integrating $f(\ell_1)$ with respect to $\ell_1$, we obtain the distribution function of $\ell_1$ as
\begin{align*}
\Pr(\ell_1<x)= \frac{\Gamma_n(n+1/2)(\frac{x}{2})^{nm/2}}{\Gamma_n\{(n+m+1)/2\}|\Sigma|^{n/2}}
\sum_{k=0}^{\infty}\sum_{\kappa \in P^k_n} \frac{(m/2)_\kappa C_\kappa(-\frac{1}{2}x\Sigma^{-1})C_\kappa(I_n)}{\{(n+m+1)/2\}_\kappa k!C_\kappa( I_m)}.
\end{align*}
The zonal polynomials $C_{\kappa}(I_m)$ are expressed, for the length of partition $p>0$, as
\begin{align} \label{eq:CIm}
C_\kappa(I_m)=\frac{2^{2k}k!(m/2)_\kappa \prod_{i<j}^{p}(2\kappa_i-2\kappa_j-i+j)}{\prod_{i=1}^{p}(2\kappa_i+p-i)!}.
\end{align}
Using the heterogeneous hypergeometric function, we obtain the distribution function of $\ell_1$ as
\begin{align*}
\Pr(\ell_1<x)&= \frac{\Gamma_n\{(n+1)/2\}(\frac{x}{2})^{nm/2}}{\Gamma_n\{(n+m+1)/2\}|\Sigma|^{n/2}}
\sum_{k=0}^{\infty}\sum_{\kappa \in P^k_n} \frac{(n/2)_\kappa C_\kappa(-\frac{1}{2}x\Sigma^{-1})}{\{(n+m+1)/2\}_\kappa k!}\\
&=\frac{\Gamma_n\{(n+1)/2\}(\frac{x}{2})^{nm/2}}{\Gamma_n\{(n+m+1)/2\}|\Sigma|^{n/2}}
%\sum_{k=0}^{\infty}\sum_{\kappa \vdash k\atop{l(\kappa)=n}}^{}
\; {_1F_1}^{(m,n)}\left(\frac{m}{2};\frac{n+m+1}{2};-\frac{1}{2}x\Sigma^{-1},I_n\right ),
\end{align*}
where we note that $(m/2)_\kappa / C_\kappa(I_m) = (n/2)_\kappa / C_\kappa(I_n)$ from \eqref{eq:CIm}.
\end{proof}
%%%%%%%%%%%%%%%%%%%%%%%%%%%%%%%%%%%%%%%%%%%%%%%%%%%
\begin{corollary}
Let $W\sim \mathcal{W}_m(n, I_m)$, with $m>n$. Then the distribution function of the largest eigenvalue $\ell_1$ of $W$ is given as
\begin{align}
\Pr(\ell_1<x)
\label{eq:eq25}
&=\frac{\Gamma_n\{(n+1)/2\}(\frac{x}{2})^{nm/2}}{\Gamma_n\{(n+m+1)/2\}}\mathrm{exp}\left(-\frac{nx}{2}\right)
{{}_1F_1}\left(\frac{n+1}{2};\frac{n+m+1}{2};\frac{x}{2}I_n\right ).
\end{align}
\end{corollary}
\noindent
\begin{proof}
From Theorem~2, we have
\begin{align}
\label{eq:eq24}
\Pr(\ell_1<x)=
\frac{\Gamma_n\{(n+1)/2\}(\frac{x}{2})^{nm/2}}{\Gamma_n\{(n+m+1)/2\}}
{{}_1F_1}\left(\frac{m}{2};\frac{n+m+1}{2};-\frac{x}{2}I_n\right ).
\end{align}
The distribution function (\ref{eq:eq24}) is translated to the series of positive terms using the Kummer relation as 
${{}_1F_1}(a;c;-X)=\mathrm{etr}(-X){{}_1F_1}(c-a;c;X)$.
\end{proof}

\begin{figure}[H]
 \begin{center}
     \includegraphics[width=7cm]{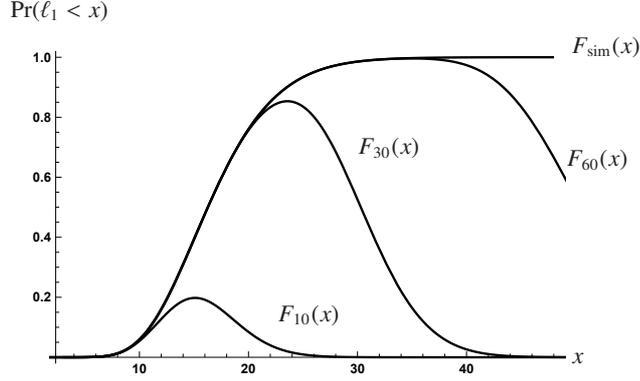}
        \rlap{\raisebox{4.5ex}{\kern-11.0em{\small $F_{10}(x)$}}}%
         \rlap{\raisebox{18.5ex}{\kern-8.0em{\small $F_{30}(x)$}}}%
         \rlap{\raisebox{17.5ex}{\kern-0.20em{\small $F_{60}(x)$}}}%
         \rlap{\raisebox{27.0ex}{\kern-0.0em{\small $F_{\mathrm{sim}}(x)$}}}%
     \rlap{\raisebox{1.0ex}{\kern-0.0em{\small $x$}}}%
\rlap{\raisebox{30.0ex}{\kern-21.0em{\small $\Pr(\ell_1 < x) $}}}%
        \caption{Comparison of $F_{\mathrm{sim}}(x)$ and $F_{K}(x)$ for $K=10,30,60$, where $m=10$, $n=3$, $\Sigma=I_{10}$.} 
        \label{fig:eigenmax}
                \end{center}
\end{figure}

%We now ready to calculate the distribution functions of (\ref{eq:eq25})．
The distribution function (\ref{eq:eq25}) is an infinite series. 
The truncated distribution function of (\ref{eq:eq25}) up to the $K$th degree is denoted by 
\begin{align}
\label{eq:eq26}
F_K(x)
&=\frac{\Gamma_n\{(n+1)/2\}(\frac{x}{2})^{nm/2}}{\Gamma_n\{(n+m+1)/2\}}\mathrm{exp}\left(-\frac{nx}{2}\right)
\sum_{k=0}^{K}\sum_{\kappa \in P^k_n}\frac{\{(n+1)/2\}_\kappa}{\{(n+m+1)/2\}_\kappa}\frac{C_\kappa(\frac{x}{2}I_n)}{k!}.
\end{align}
The empirical distributions based on $10^6$-trial Monte Carlo simulations are denoted by $F_{\mathrm{sim}}(x)$. 
Fig. \ref{fig:eigenmax} shows the comparison of $F_{\mathrm{sim}}(x)$ and $F_{K}(x)$ for $K=10,30,60$. 
If $K=60$, then the truncated series $F_{K}(x)$ reaches that near $x=40$. 
Table~\ref{table:pp-l1} shows the comparison of percentile points between $F^{-1}_{\text{sim}}$, $F^{-1}_{K}$, and $F^{-1}$, where the exact $F^{-1}$ is 
calculated by the method of \cite{Chiani2014}. All percentile points have the same precision.
In the case of $n=3$, we need about $60$ and $90$ terms in the hypergeometric series for $m=10$ and $m=50$, respectively. 
However, more terms in the hypergeometric series are needed and much longer calculation time is required compared with the case for $n=2$.
We observe that the hypergeometric series in (\ref{eq:eq26}) converges slowly when the sample of size $n$ increases. 
%%%%%%%%画像
%%%%%%%%%%%%%%%%%%%%%%%%%%%%
\begin{table}[H]
\caption{Comparison of percentile points of $\ell_1$ of $W \sim \mathcal{W}_{m}{(3, I_m})$ between $F^{-1}_{\text{sim}}$, $F^{-1}_{K}$, and $F^{-1}$ for dimension $m=10$ and $m=50$.} \label{table:pp-l1}
\begin{center}
\begin{tabular}{c}
    \begin{minipage}[c]{0.4\hsize}
      \begin{center}
%\tbl{Comparison of acoustic for frequencies for piston-cylinder problem.}
       \captionsetup{labelformat=empty,labelsep=none}
          \subcaption{$m=10$}
{\begin{tabular}{@{}cccc@{}} \toprule
$\alpha$&${{F}^{\small{-1}}_{\mathrm{sim}}}$ &$F^{-1}_{60}$&$F^{-1}$  \\
0.01	 &~7.75& 	7.75&7.75\\
0.05	 &~9.74& 	9.74&9.74\\
0.50	 &~16.2& 	16.2&16.2\\
0.95	 &~25.9& 	25.9&25.9\\
0.99  &~31.1& 	31.1&31.1\\
\noalign{\smallskip}\hline
\end{tabular}}
 \end{center}
  \end{minipage}
  \begin{minipage}[c]{0.4\hsize}
          \begin{center}
        \captionsetup{labelformat=empty,labelsep=none}
         \subcaption{$m=50$}
{\begin{tabular}{@{}cccc@{}} \toprule
$\alpha$&${{F}^{\small{-1}}_{\mathrm{sim}}}$&${{F}^{\small{-1}}_{90}}$  &$F^{-1}$ \\
0.01	&46.2&~46.2&46.2 	\\
0.05	 &50.9&~50.9& 50.9\\
0.50	 &64.4&~64.4&64.4 \\
0.95	 &81.3&~81.4&81.4 \\
0.99  &89.6&~89.7& 89.7\\
\noalign{\smallskip}\hline
\end{tabular}}
        \end{center}
    \end{minipage}
  \end {tabular}
  \end{center}
\end{table}
We consider the non-null case of $W\sim \mathcal{W}_1(2,\Sigma)$, where $\Sigma \neq I_2$.
The density function of eigenvalue $\ell_1$ of $W$ is given as  
\begin{align}
\label{eq:eq28}
f(\ell_1)&= \frac{1}{2 \sqrt{|\Sigma|}}\exp \biggl(-\frac{1}{2\lambda_1}\ell_1\biggl){_1F_1}\biggl(\frac{1}{2};1;a\biggl),
\end{align}
where $\Sigma=\mathrm{diag}(\lambda_1, \lambda_2)$ and $a=-\frac{1}{2} \ell_1 (1/\lambda_2 - 1/\lambda_1)$.
The derivation of (\ref{eq:eq28}) is presented in the appendix.
 Fig. \ref{fig:L1-m=2} (a) shows the comparison of (\ref{eq:eq28}) and $10^6$ Monte Carlo simulation results 
 under $\Sigma=\mathrm{diag}(5,2)$.
  Fig. \ref{fig:L1-m=2} (b) shows the line for $F_{100}(x)$ and the dot plots for numerical integration of (\ref{eq:eq28}).
  They have almost the same precision.
  %%%%%%%%%%%%%%%%%%%%%%%%%%
  %%%%%画像
\begin{figure}[H]
\begin{center}
\begin{tabular}{cc}
\includegraphics[width=7cm]{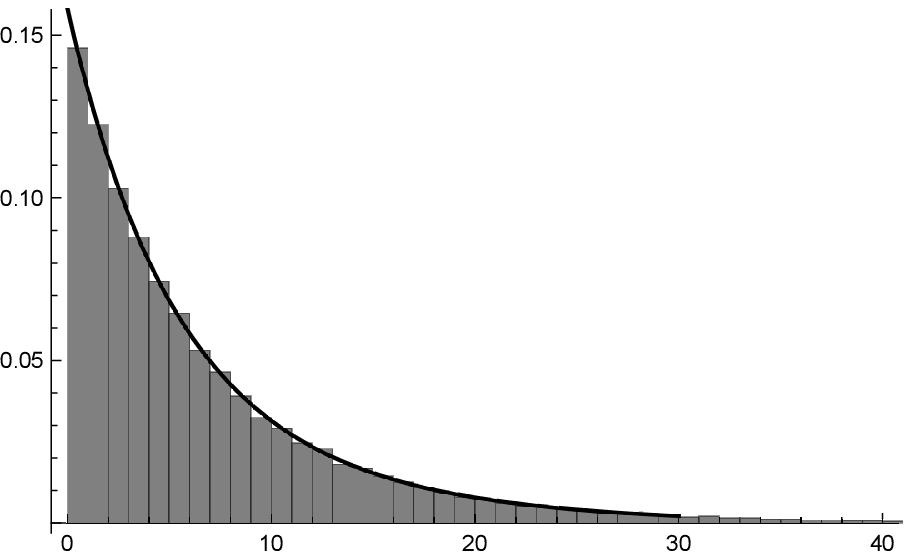} & 
\includegraphics[width=7cm]{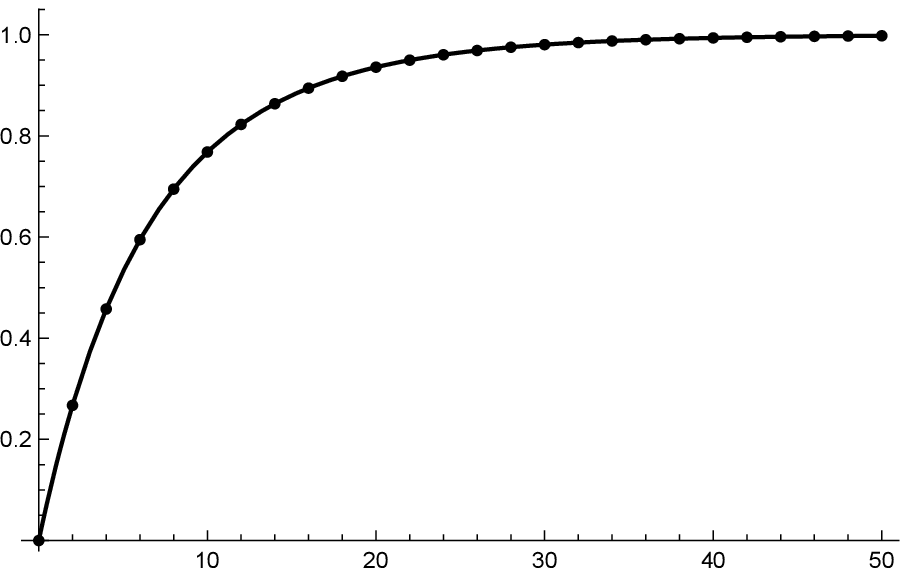}
\rlap{\raisebox{1.0ex}{\kern-21.0em{\small $x$}}}%
\rlap{\raisebox{1.0ex}{\kern-0.0em{\small $x$}}}%
\rlap{\raisebox{30.0ex}{\kern-21.0em{\small $\Pr(\ell_1 < x) $}}}%
\rlap{\raisebox{30.0ex}{\kern-41.0em{\small$f(\ell_1)$}}}%
\\
(a) $f(\ell_1)$ in \eqref{eq:eq28} and Monte Carlo Sim. for $\ell_1$ & (b) $F_{100}(x)$ and numerical integrate of (\ref{eq:eq28})
\end{tabular}
\caption{Comparison of the theoretical distribution of $\ell_1$ of $W \sim \mathcal{W}_2(1, \mathrm{diag}(5,2))$ and Monte Carlo simulation.}
\label{fig:L1-m=2}
\end{center}
\end{figure}
%%%%%%%%%%%%%%%%%%%%%%%%%%%%%%%%%%%%%%%%%%%%%%%%%%%%%%%%%%%%%%%%%
\section{Stiefel manifold over a real finite-dimensional division algebra}
\label{sec:04}
%%%%%%%%%%%%%%%%%%%%%
\if0
The eigenvalues of real and complex Wishart matrices have been used in multivariate analysis for various applications. Recently, this class has been generalized to any $\beta>0$ that includes classical results for real, complex, and quaternion cases．
Koev and Edelman~(2013) introduced the beta-Wishart ensemble and derived the exact distributions of the extreme eigenvalues for a non-singular beta-Wishart matrix. The essence of this approach is the use of Jack polynomials instead of zonal polynomials. On the other hand, D\'iaz-Garc\'ia~(2011) expanded a Wishart ensemble on a normed division algebra. D\'iaz-Garc\'ia~(2013) later derived some integral formulas for Jack polynomials over the Stiefel manifold on normed division algebra. We also define the heterogeneous hypergeometric functions of parameter $\beta>0$ and present some useful formulas for them. 
\fi
%%%%%%%%%%%%%%%%%%
%In this section, we define the heterogeneous hypergeometric functions of parameter $\beta>0$ in terms of Jack polynomials on division algebra. 
%Gross and Richards~(1987) studied the hypergeometric functions on division algebra. 
Let $\mathbb{F}_\beta$ denote a real finite-dimensional division algebra such that
$\mathbb{F}_1 = \mathbb{R}$, $\mathbb{F}_2 =\mathbb{C}$, and $\mathbb{F}_4 = \mathbb{H}$ for $\beta=1, 2, 4$, where
$\mathbb{R}$ and $\mathbb{C}$ are the fields of real and complex numbers, respectively, and $\mathbb{H}$ is the quaternion division algebra over $\mathbb{R}$.
We restrict the parameter $\beta$ to values of $\beta = 1, 2$, and $4$, and denote $\mathbb{F}_\beta^{m \times n}$ 
by the set of all $m \times n$ matrices over $\mathbb{F}_\beta$, where $m \ge n$.
The conjugate transpose of $X \in \mathbb{F}_\beta^{m \times n}$ is written by $X^\ast=\overline{X^\top}$ and we say that $X$ is Hermitian if $X^\ast = X$.  The set of all Hermitian matrices is denoted by 
$S^{\beta}(m)=\{ X \in \mathbb{F}^{m\times n}_\beta \mid X^\ast=X\}$. 
The eigenvalues of a Hermitian matrix are all real.  If the eigenvalues of $X \in S^\beta(m)$ are all positive, 
then we say that it is positive definite and write $X > 0$.
The exterior product $(dX)$ for $X \in \mathbb{F}_\beta^{m \times n}$ was defined in Mathai~\cite{Mathai1997} and D\'iaz-Garc\'ia and Guti\'errez-J\'aimez~\cite{Garcia2011}. 
In a similar manner (\ref{stiefel:def}), we define the Stiefel manifold and the unitary group over $\mathbb{F}_\beta$ as
\begin{align*}
V^{\beta}_{n,m} = \{ H_1\in \mathbb{F}_\beta^{m\times n}\mid H_1^{\ast}H_1=I_n\},
\quad  U^{\beta}_m = V^{\beta}_{m,m}=\{ H \in \mathbb{F}_\beta^{m\times m} \mid H^{\ast}H = H H^{\ast} = I_m\} ,
\end{align*}
respectively.
If $\beta \in\{1, 2, 4\}$, then $U^{\beta}_m$ are the real orthogonal group, unitary group, and symplectic group, respectively. 
The $\beta$-multivariate gamma function for $c \in \mathbb{F}_\beta$, $\Gamma_m^{\beta}(c)$, is defined by
\begin{align*}
\Gamma_m^{\beta}(c) &= \int_{X > 0} |X|^{c - (m-1)\beta/2 } \mathrm{etr}(-X) (d X)
= \pi^{\frac{m(m-1)\beta}{4}}\prod_{i=1}^{m}\Gamma\bigg\{c-\frac{(i-1)\beta}{2}\bigg\}, 
\end{align*}%ガンマ関数定義
where $\mathrm{Re}(c)>(m-1)\beta/2$.
We define $(H_1^{\ast}dH_1)$ and $\mathrm{Vol}(V^{\beta}_{n,m})$ by
\begin{align*}
(H_1^{\ast}dH_1)=\bigwedge_{i=1}^{n}\bigwedge_{j=i+1}^{m}h_j^{\ast}dh_i,
\quad 
\mathrm{Vol}(V^{\beta}_{n,m})=\int_{H_1\in V^{\beta}_{n,m}}(H_1^{\ast}dH_1)=\frac{2^n\pi^{mn\beta/2}}{\Gamma^{\beta}_n(m\beta/2)} ,
\end{align*}
respectively, where 
$H_1\in V^{\beta}_{n,m}$ and $H = (H_1 \mid H_2)=(h_1, \dots, h_n \mid h_{n+1}, \dots, h_m) \in U_m^\beta$.
Another differential form $(dH_1)$ defined by
\begin{align*}
(dH_1)=\frac{(H_1^{\ast}dH_1)}{\mathrm{Vol}(V^{\beta}_{n,m})}=\frac{\Gamma^{\beta}_n(m\beta/2)}{2^n\pi^{mn\beta/2}}(H_1^{\ast}dH_1)
\end{align*}
is normalized such as 
$
\int_{H_1\in V^{\beta}_{n,m}}(dH_1)=1.
$
For a partition $\kappa$, the $\beta$-generalized Pochhammer symbol of parameter $a>0$ is defined as
\begin{align*}
(a)_\kappa^{\beta}=\prod_{i=1}^{m}\biggl(a-\frac{i-1}{2}\beta\biggl)_{\kappa_i}.
\end{align*}

The Jack polynomial $C_\kappa^{\beta}(X)$ is a symmetric polynomial in $x_1,\dots,x_m$; these are eigenvalues of $X$.
 See Stanley~\cite{Stanley1989} and Koev and Edelman~\cite{Edelman2006} for the relevant detailed properties.
If $\beta=1,2$, then Jack polynomials are referred to as zonal polynomials and Shur polynomials, respectively.
Li and Xue~\cite{Li2009} proposed zonal polynomials and hypergeometric functions of quaternion matrix arguments for $\beta=4$.
For $A\in S^{\beta}(m)$ and $B\in S^{\beta}(n)$, the heterogeneous hypergeometric functions of parameter $\beta$ are defined as 
\begin{align}
\label{def:2}
_pF_q^{(\beta;m,n)}({\mbox{\boldmath$\alpha$}};\mbox{\boldmath$\beta$};A,B)=\sum_{k=0}^{\infty}\sum_{\kappa \in P^k_n}\frac{(\alpha_1)^{\beta}_\kappa \cdots(\alpha_p)^{\beta}_\kappa}{(\beta_1)^{\beta}_\kappa \cdots(\beta_q)^{\beta}_\kappa}\frac{C^{\beta}_\kappa(A)C^{\beta}_\kappa(B)}{k!C^{\beta}_\kappa(I_m)}, 
\end{align}
where $m\geq n$. 
For $B\in S^{\beta}(n)$, the functions $_pF_q^{(\beta;n)}({\mbox{\boldmath$\alpha$}},\mbox{\boldmath$\beta$};B)$ are also defined as $_pF_q^{(\beta;m,n)}({\mbox{\boldmath$\alpha$}},\mbox{\boldmath$\beta$};I_m,B)$. 
If $\beta=1$, we write the functions $_pF_q({\mbox{\boldmath$\alpha$}},\mbox{\boldmath$\beta$};B)$ by $_pF_q^{(1;n)}({\mbox{\boldmath$\alpha$}},\mbox{\boldmath$\beta$};B)$.
Gross and Richards~\cite{Gross1987} provided some properties of the integral formula for hypergeometric functions on a division algebra. 
The following two equations are needed for the derivation of Theorems~\ref{th:s04-01} and \ref{th:s04-02}, respectively. 
We utilize the complexification $\mathbb{C}S^{\beta}(m)=S^{\beta}(m)+iS^{\beta}(m)$. 
That is, $\mathbb{C}S^{\beta}(m)$ consists of all matrices $Z$ of the form $Z=X+iY$, where $X, Y\in S^{\beta}(m)$, and $i^2=-1$.
We also write $\mathrm{Re}(Z) = X$ and $\mathrm{Im}(Z) = Y$ for $Z=X+iY$.

For $Y\in \mathbb{C}S^{\beta}(n)$, $t=(n-1)\beta/2$, $\mathrm{Re}(a)>(n-1)\beta/2$ and $\mathrm{Re} (b)>(n-1)\beta/2$, we have
\begin{align}
\label{eqn:15}
\int_{{0}<U<{I_n}}|U|^{a-t-1}|I_n-U|^{b-t-1}C^{\beta}_{\kappa}(UY)(dU)
&=\frac{\Gamma^{\beta}_n(a,\kappa)
\Gamma^{\beta}_n(b)
}{\Gamma^{\beta}_n(a+b,\kappa)
}C^{\beta}_{\kappa}(Y).
\end{align}
For $Y\in S^{\beta}(m)$ and $Z\in \mathbb{C}S^{\beta}(n)$ with $\mathrm{Re}(Z)>0$, we have 
\begin{align}
\label{eqn:16}
\int_{U>0}\mathrm{etr}(-UZ)|U|^{a-t-1}C^{\beta}_\kappa(UY)(dU)&=(a)^{\beta}_\kappa\Gamma^{\beta}_n(a)|Z|^{-a}C^{\beta}_\kappa(YZ^{-1}). 
\end{align}
The integral representation for the functions ${_{1}F_{1}}^{(\beta;m,n)}$ and ${_{2}F_{1}}^{(\beta;m,n)}$ in (\ref{def:2}) are given in Theorem~\ref{th:s04-01}.
%%%%%%%%%%%%%%%%%%%%%%%%%%%%%%%%%%%%%%%%%%%%%%%%%%%%%%%%%%%%%%%%%．
\begin{theorem} \label{th:s04-01}
For $X\in S^{\beta}(m)$ and $Y\in \mathbb{C}S^{\beta}(n)$, the function ${_{1}F_{1}}^{(\beta;m,n)}$ is represented by the following integral representation:
\begin{align*}
&{_{1}F_{1}}^{(\beta;m,n)}(a;c;X,Y)=\frac{\Gamma^{\beta}_n(c)}{\Gamma^{\beta}_n(a)\Gamma^{\beta}_n(c-a)}\int_{{0}<U<{I_n}}{_0F_0}^{(\beta;m,n)}(X,UY)|U|^{a-t-1}|I_n-U|^{c-a-t-1}(dU),
\end{align*}
where $\mathrm{Re}(c)>\mathrm{Re}(a)+(n-1)\beta/2>(n-1)\beta$. 

For arbitrary $a_1$, the function ${_{2}F_{1}}^{(\beta;m,n)}$ also has the integral representation 
\begin{align*}
&{_{2}F_{1}}^{(\beta;m,n)}(a_1,a;c;X,Y)=\frac{\Gamma^{\beta}_n(c)}{\Gamma^{\beta}_n(a)\Gamma^{\beta}_n(c-a)}\int_{{0}<U<{I_n}}{_1F_0}^{(\beta;m,n)}(a_1;X,UY)|U|^{a-t-1}|I_n-U|^{c-a-t-1}(dU),  \nonumber
\end{align*}
where $\mathrm{Re}(c)>\mathrm{Re}(a)+(n-1)\beta/2>(n-1)\beta$ and $||Y||<1$, where $||Y||$ is the maximum of the absolute values of the eigenvalues of $Y$. 
\end{theorem}
\noindent
\begin{proof}
The desired result is obtained by expanding ${_0F_0}^{(\beta;m,n)}$ and ${_1F_0}^{(\beta;m,n)}$ in integrand and integrating term by term using identities (\ref{eqn:15}), respectively.
\end{proof}

%%%%%%%%%%%%%%%%%%%%%%%%%%%%%%%%%%%%%%%%%%%%%%%%%%%%%%%%%%%%%%%%%
\begin{theorem} \label{th:s04-02}
For $X\in S^{\beta}(m)$, $Y\in S^{\beta}(n)$, and $Z \in \mathcal{U}^{\beta}(n)$ with $\mathrm{Re(Z)>0}$, we have
\begin{align*}
\int_{U>0}\mathrm{etr}(-UZ)|U|^{a-t-1}{_{p}F_q}^{(\beta;m,n)}({\mbox{\boldmath$\alpha$}};\mbox{\boldmath$\beta$};X,UY)(dU) =\Gamma^{\scalebox{0.5}{$(\beta)$}}_n(a)|Z|^{-a}{_{p+1}F_q}^{(\beta;m,n)}({\mbox{\boldmath$\alpha$}},a;\mbox{\boldmath$\beta$};X,YZ^{-1}), 
\end{align*}
where $p<q$, $\mathrm{Re}(a)>(n-1)\beta/2;$ or $p=q$, $\mathrm{Re}(a)>(n-1)\beta/2$, and $||Z||^{-1}<1$. 
\end{theorem}
\noindent
\begin{proof}
The result is immediately obtained by expanding ${_{p}F_q}^{(\beta;m,n)}$ in the integrand and integrating term by term using identities (\ref{eqn:16}). 
\end{proof}

The integral formula over the Stiefel manifold in Section~\ref{th:s04-02} can be extended to the general case of division algebra. Concerning the equations (\ref{splitting}) and (\ref{eq:eq35}), these formulas can 
be extended easily to the following equations. For $A\in S^{\beta}(m)$ and $B\in S^{\beta}(n)$, we have 
\begin{align}
\label{splitting:beta}
  \int_{H_1\in V^{\beta}_{n,m}}^{}C^{\beta}_{\kappa}(AH_1BH_1^\ast)(dH_1)&=\frac{C^{\beta}_{\kappa}(A)C^{\beta}_{\kappa}(B)}{C^{\beta}_{\kappa}(I_m)}
\end{align}
and
\begin{align} \label{intCAHBH}
\int_{H_1\in V^\beta_{n,m}}C^\beta_\kappa(H_1BH_1^\ast)(dH_1)
&=\int_{H\in U^\beta_m}C^\beta_\kappa(HB_1H^\ast)(dH), 
\end{align}
where 
$B_1=\left( 
\begin{array}{cc}
B & O\\
O & O
\end{array}
\right).$
D\'iaz-Garc\'ia~\cite{Garcia2013} showed that the denominator on the right side of (\ref{splitting:beta}) was evaluated as 
 $C^{\beta}_{\kappa}(I_r)$ instead of $C^{\beta}_{\kappa}(I_m)$, where $r=\mathrm{rank}(B)$.
If $A=I_m$ in (\ref{splitting:beta}), then (\ref{splitting:beta})  and (\ref{intCAHBH}) imply 
the following well-known property of Jack polynomials:
\begin{align*}
%\int_{H_1\in V_{n,m}}C^\beta_\kappa(H_1BH_1')(dH_1)
%&=&
C^\beta_\kappa(B_1) = \int_{H\in U^\beta_m}C^\beta_\kappa(HB_1H^\ast)(dH).
\end{align*}

%%%%%%%%%%%%%%%%%%%%%%%%%%%%%%%%%%%%%%%%%%%%%%%%%%%%%%%%%%%%%%%%%
\begin{theorem}
For $A\in S^{\beta}(m)$, $B\in S^{\beta}(n)$, and $H_1\in V^{\beta}_{n,m}$, 
where $H=(H_1,H_2)\in U^{\beta}(m)$ and $\alpha_1$ is a non-negative integer, then we have
\begin{align*}
  {_pF_q}^{(\beta;m,n)}(\mbox{\boldmath$\alpha$};\mbox{\boldmath$\beta$};A,B)&=
 \displaystyle \int_{H_1\in V^{\beta}_{n,m}}{_pF_q}^{(\beta;m)}(\mbox{\boldmath$\alpha$};\mbox{\boldmath$\beta$};AH_1BH_1^\ast)(dH_1),
 \quad{_0F_0}^{(\beta;m,n)}\left(A,B\right)=
 \displaystyle \int_{H_1\in V^{\beta}_{n,m}}\mathrm{etr} \left(AH_1BH_1^\ast\right)(dH_1),\\
 {_1F_0}^{(\beta;m,n)}(\alpha_1;A,B)&=\int_{H_1\in V^{\beta}_{n,m}}|I_m-AH_1BH_1^\ast|^{-\alpha_1}(dH_1), 
 \quad{_0F_0}^{(\beta;m,n)}\left(I_m+A,B\right)=\mathrm{etr}(B)~ {_0F_0}^{(\beta;m,n)}\left(A,B\right).
    \end{align*}
    \end{theorem}
   \noindent
    \begin{proof}
    The derivation is the same as that for the real case in Section~\ref{sec:02}.
    \end{proof}

\section{Singular beta-Wishart matrix}
\label{sec:05}
In this section, we define the singular  beta-Wishart matrix on a real finite division algebra. This matrix covers the singular real, complex, and quaternion Wishart matrices. We derive the density function of the singular beta-Wishart distributions and some distributions of eigenvalues. The singular beta-Wishart distributions are denoted by $\mathcal{W}^{\beta}_m(n,\Sigma)$. 

Let an $m\times n$ beta-Gaussian random matrix $X$ be distributed as $X\sim$ $\mathcal{N}^{\beta}_{m,n}(M,\Sigma \otimes \Theta)$, where $M$ is an $m \times n$ mean matrix and $\Sigma$ and $\Theta$ are $m \times m$ and $n \times n$ positive definite matrices, respectively, 
The density functions of $X$ are given as 
\begin{align*}
 \frac{1}{(2\pi \beta^{-1})^{mn\beta/2}|\Sigma|^{\beta n/2}|\Theta|^{\beta m/2}}\mathrm{exp}\biggl(-\frac{\beta}{2}\mathrm{tr}\Sigma^{-1}(X-M)\Theta^{-1}(X-M)^{\ast}\biggl).
\end{align*}
Let $X\sim \mathcal{N}^{\beta}_{m,n}({O},\Sigma \otimes I_n)$; that is, $M={O}$ and $\Theta=I_n$. 
Then the $m\times m$ singular beta-Wishart matrix is defined as $W=XX^{\ast}$, where $m> n$.
Real~($\beta=1$) and complex~($\beta=2$) singular Wishart distributions were obtained by Uhlig~\cite{Uhlig1994},
Srivastava~\cite{Srivastava2003}, and Ratnarajah and Vaillancourt~\cite{Ratnarajah2005a}. D\'iaz-Garc\'ia and Guti\'errez-S\'anchez~\cite{Garcia2013andGutierrez-Sanchez} derived some useful Jacobians of the transformation for singular matrices. 
Let $W\sim \mathcal{W}^{\beta}_m(n,\Sigma)$. From the Jacobian of the transformation of Corollary~1 in D\'iaz-Garc\'ia and Guti\'errez-S\'anchez~\cite{Garcia2013andGutierrez-Sanchez},  the density function of $W$ is given as 
\begin{align*}
 f(W)=\frac{\pi^{n\beta(n-m)/2}(\mathrm{det}\Sigma)^{-\beta n/2}}{(2\beta^{-1})^{\beta mn/2}\Gamma^{\beta}_n(n\beta/2)}\mathrm{etr}\biggl(-\frac{\beta}{2}\Sigma^{-1}W\biggl)(\mathrm{det}L_1)^{\beta(n-m+1)/2-1},
 \end{align*}
 where $H_1\in V_{n,m}^{\beta}$, $L_1=\mathrm{diag}(\ell_1,\dots,\ell_n)$ and  $W=H_1L_1H_1^{\ast}$.
 %%%%%%%%%%%%
Srivastava~\cite{Srivastava2003} gave the joint density as 
 \begin{align}
\label{joint:beta}
 f(\ell_1,\dots,\ell_n)\propto \int_{H_1\in V^{1}_{n,m}}\mathrm{etr}\biggl(-\frac{1}{2}\Sigma^{-1}H_1L_1H_1^\top\biggl)(dH_1). 
 \end{align}
 However, the right hand side of (\ref{joint:beta}) is not represented in terms of heterogeneous hypergeometric functions.
 On the other hand, Ratnarajah and Vaillancourt~\cite{Ratnarajah2005a} gave the joint density of eigenvalues of a singular complex Wishart matrix as 
  \begin{align}
  \label{hyper0f0}
 f(\ell_1,\dots,\ell_n)\propto \int_{H_1\in V^{2}_{n,m}}\mathrm{etr} \left(-\Sigma H_1L_1H_1^\ast\right)(dH_1)
 ={_0F_0}^{(2;m,n)}\left(-\Sigma,L_1 \right).
 \end{align}
 However, they did not prove the equation given above. 
 From (\ref{splitting:beta}), the equation (\ref{hyper0f0}) is easily proved. 
  Theorem~6 presents the joint density of eigenvalues of a singular beta-Wishart matrix as follows.
  The joint density of eigenvalues (\ref{eq:eq13}) can be expanded on a division algebra. 
  %%%%%%%%%%%%%%%%%%%%%%%%%%%%%%
  \begin{theorem}
 Let $W\sim \mathcal{W}^{\beta}_m(n,\Sigma)$; then the joint density of eigenvalues $\ell_1,\dots,\ell_n$ of $W$ is given as 
 \begin{align}
   \label{jointeigen01}
 \displaystyle {f(\ell_1,\dots,\ell_n)=\frac{(2\beta^{-1})^{-\beta nm/2}\pi^{n^2\beta/2+r}}{|\Sigma|^{\beta n/2}\Gamma^{{\beta}}_n(\frac{n\beta}{2})\Gamma^{{\beta}}_n(\frac{m\beta}{2})}(\mathrm{det}L_1)^{\beta(m-n+1)/2-1}} \displaystyle \prod_{i<j}^{n}(\ell_i-\ell_j)^\beta{_0F_0}^{(\beta;m,n)}\biggl(-\frac{\beta}{2}\Sigma^{-1},L_1\biggl),
\end{align}
where $m>n$ and \begin{align*}r&=
\begin{cases}
\quad 0,\quad  \quad~~~~\text{\ $\beta=1,$}\\
\quad -\beta n/2,~~~\text{$\beta=2,4$}.\\
\end{cases}
\end{align*}
 \end{theorem}
 \noindent
\begin{proof}
The Jacobian of the transformation $W=H_1L_1H_1^{\ast}$ given in D\'iaz-Garc\'ia and Guti\'errez-S\'anchez~\cite{Garcia2013andGutierrez-Sanchez} is 
\begin{align}
\label{jacobian01}
(dW)=2^{-n}\pi^{r}\prod_{i=1}^{n}\ell_i^{\beta(m-n)}\prod_{i<j}^{n}(\ell_i-\ell_j)^{\beta}(dL)\wedge(H_1^{\ast}dH_1). 
\end{align}
Using identities (\ref{jacobian01}) for the density function of a singular beta-Wishart matrix $f(W)$ and integrating with respect to $H_1$ over the Stiefel manifold $V_{n,m}^{\beta}$, we have 
 \begin{align*}
\displaystyle {f(\ell_1,\dots,\ell_n)=\frac{(2\beta^{-1})^{-\beta nm/2}\pi^{n^2\beta/2+r}}{|\Sigma|^{\beta n/2}\Gamma^{\beta}_n(\frac{n\beta}{2})\Gamma^{\beta}_n(\frac{m\beta}{2})}(\mathrm{det}L_1)^{\beta(m-n+1)/2-1}} \displaystyle \prod_{i<j}^{n}(\ell_i-\ell_j)^\beta \int_{H_1\in V^{\beta}_{n,m}}\mathrm{etr}\biggl(-\frac{\beta}{2}\Sigma^{-1}H_1L_1H_1^{\ast}\biggl)(dH_1)  .
 \end{align*}
 From Theorem~5, we have the desired result. 
\end{proof}

 To derive the exact distributions of the largest eigenvalue of a singular beta-Wishart matrix, we extend Lemma~2 to the case of division algebra.  
 %%%%%%%%%%%%%
   \begin{lemma}
   Let $L=\mathrm{diag}(\ell_1,\dots ,\ell_n)$ and let the length of $\kappa$ be equal to or less than $n$; 
then the following equation holds;
    \begin{align}
   \label{sugiyama:lem03}
\int_{1>\ell_1>\ell_2>\cdots \ell_n>0}|L|^{a-t-1}|I_n-L|^{b-t-1}C^{{\beta}}_\kappa(L)\prod_{i<j}^{n}(\ell_i-\ell_j)^{\beta}\prod_{i=1}^{n}d\ell_i
=\frac{\Gamma^{{\beta}}_n(n\beta/2)\Gamma^{{\beta}}_n(a,\kappa)\Gamma^{{\beta}}_n(b)C^{{\beta}}_\kappa(I_n)}{\pi^{n^2\beta/2+r}\Gamma^{{\beta}}_n(a+b,\kappa)}.
\end{align}
  \end{lemma}
 \noindent
\begin{proof}
Let $Y=I_n$ and $U=HLH^{\ast}$ in (\ref{eqn:15}). From Proposition~3 in D\'iaz-Garc\'ia and Guti\'errez-J\'aimez~\cite{Garcia2011}, the differential form $(dU)$ is represented as
\begin{align*}
(dU)=2^{-n}\pi^{r}\prod_{i<j}^{n}(\ell_i-\ell_j)^{\beta}(dL)(H^{\ast}dH).
\end{align*}
Using the above differential form $(dU)$ and integrating $(H^{\ast}dH)$ with respect to $H$ over $U^{\beta}_n$, we have the desired result.
\end{proof}
%%%%%%%%%%%%%%
  \begin{theorem}
  Let $X_1=\mathrm{diag}(1,x_2,\dots,x_n)$ and $X_2=\mathrm{diag}(x_2,\dots,x_n)$ with $x_2>\cdots >x_n>0$; then the following equation holds.
  \begin{align}
  \nonumber
   \label{sugiyama:lem04}
\int_{1>x_2>\cdots x_n>0}|X_2|^{a-t-1}C^{{\beta}}_\kappa(X_1)\prod_{i=2}^{n}(1-x_i)^\beta \prod_{i<j}(x_i-x_j)^\beta \prod_{i=2}^{n}dx_i\\
=(na+k)(\Gamma^{{\beta}}_n(\beta n/2)/\pi^{n^2\beta/2+r})\frac{\Gamma^{{\beta}}_n(a,\kappa)\Gamma^{{\beta}}_n(t+1)C^\beta_\kappa(I_n)}{\Gamma^{{\beta}}_n(a+t+1,\kappa)}.
\end{align}
  \end{theorem}
 \noindent
\begin{proof}
Let $b=t+1$. With the translation of $x_i=\ell_i/\ell_1$, $i\in\{2,\ldots, n\}$, the left side of (\ref{sugiyama:lem03}) is given as 
\begin{align*}
\int_{0}^{1}\ell_1^{na+k-1}d\ell_1 \int_{1>x_2>\cdots x_n>0}|X_2|^{a-t-1}C^{{\beta}}_\kappa(X_1)\prod_{i=2}^{n}(1-x_i)^\beta \prod_{i<j}(x_i-x_j)^\beta \prod_{i=2}^{n}dx_i.
\end{align*}
 We note that $\int_{0}^{1}\ell_1^{na+k-1}d\ell_1=1/na+k$. 
\end{proof}
%%%%%%%%%%%%%
   \begin{theorem}
    Let $W\sim \mathcal{W}^{\beta}_m(n,\Sigma)$,where $m>n$. Then the distribution function of the largest eigenvalue $\ell_1$ of $W$ is given as 
   \begin{align}
     \label{betamax}
\displaystyle\mathrm{Pr}(\ell_1<x)=\frac{\Gamma^{{\beta}}_n\{(n-1)\beta/2+1\}(x\beta/2)^{\beta nm/2}}{\Gamma^\beta_n\{(n+m-1)\beta/2+1\}|\Sigma|^{n\beta/2}} 
{{}_1F_1{}}^{(\beta;m,n)}\left(\frac{m\beta}{2};\frac{(n+m-1)\beta}{2}+1;-\frac{\beta}{2}x\Sigma^{-1},I_n\right ).
\end{align}
    \end{theorem}
   \noindent
\begin{proof}
This proof is presented in the same way as that for Theorem~2 in Section~3. 
We consider the joint density of eigenvalues (\ref{jointeigen01}), the translation of $x_i=\ell_i/\ell_1$, $i\in\{2,\ldots, n\}$, and (\ref{sugiyama:lem04}) in order to integrate $x_2,\dots,x_n$ in $(\ref{jointeigen01})$. Moreover, integrating the density function $f(\ell_1)$ with respect to $\ell_1$, the density function of $\ell_1$ is obtained by
  \begin{align*}
\displaystyle\mathrm{Pr}(\ell_1<x)=\frac{\Gamma^{{\beta}}_n\{(n-1)\beta/2+1\}(x\beta/2)^{nm\beta/2}}{\Gamma^{\beta}_n\{(n+m-1\}\beta/2+1)|\Sigma|^{n\beta/2}} \sum_{k=0}^{\infty}\sum_{\kappa \in P^k_n} \frac{(m\beta/2)^\beta_\kappa \;C^{{\beta}}_\kappa(-\frac{\beta}{2}x\Sigma^{-1})\;C^{{\beta}}_\kappa(I_n)}{\{(n+m-1)\beta/2+1\}^\beta_\kappa k!C^{{\beta}}_\kappa( I_m)}.
\end{align*}
The Jack polynomials $C^{\beta}_\kappa(I_m)$ are expressed as 
\begin{align}
\label{jacknull}
C^{\beta}_\kappa(I_m)=\frac{(2\beta)^{2k}k!}{j_\kappa}\biggl(\frac{m}{2\beta}\biggl)^\beta_\kappa ,
\end{align}
where 
\begin{eqnarray*}
j_\kappa=\prod_{(i,j)\in \kappa}h_{\ast}^{\kappa}(i,j)h_{\kappa}^{\ast}(i,j)
\end{eqnarray*}
and $h_{\kappa}^{\ast}(i,j)\equiv \kappa_j'-i+2\beta(\kappa_i-j+1)$ and $h_{\ast}^{\kappa}(i,j)\equiv \kappa_j'-i+1+2\beta(\kappa_i-j)$ are the upper and lower hook lengths at $(i,j)\in \kappa$, respectively. See Koev and Edelman~\cite{Edelman2006} for details. 
Using the identities $(\ref{jacknull})$, we obtain (\ref{betamax}).
 \end{proof}

Furthermore, the following corollary holds by applied the $\beta$-Kummer relation to ${{}_1F_1}^{(\beta;m)}$ in (\ref{betamax}).

\begin{corollary}
Under the same condition of Theorem~8, the distribution function of $\ell_1$ is given by
\begin{align}
\displaystyle\mathrm{Pr}(\ell_1<x) =& \frac{
\Gamma^{{\beta}}_n\{(n-1)\beta/2+1\}(x\beta/2)^{ nm\beta/2}}{\Gamma^\beta_n\{(n+m-1)\beta/2+1\}|\Sigma|^{n\beta/2}} 
 \mathrm{etr}\biggl(-\frac{\beta x}{2}\Sigma^{-1}\biggl){
{}_1F_1{}}^{(\beta;m)}\left(\frac{(m-1)\beta}{2}+1;\frac{(n+m-1)\beta}{2}+1;\frac{\beta}{2}x\Sigma^{-1}\right ).
\label{eq:-cdf-l1-2}
\end{align}
\end{corollary}
\noindent
\begin{proof}
If the length of a partition $\kappa$ is $m$, then we have $(\beta n / 2)^\beta_\kappa=0$, where $m > n$.
Therefore  the function ${{}_1F_1{}}^{(\beta;m,n)}$ in (\ref{betamax}) is rewritten by
\begin{align}
\nonumber 
& {{}_1F_1{}}^{(\beta;m,n)}\left(\frac{m\beta}{2};\frac{(n+m-1)\beta}{2}+1;-\frac{\beta}{2}x\Sigma^{-1},I_n\right )
=\sum_{k=0}^{\infty}\sum_{\kappa \in P^k_n} \frac{(m\beta/2)^\beta_\kappa \;C^{{\beta}}_\kappa(-\frac{\beta}{2}x\Sigma^{-1})C^{{\beta}}_\kappa(I_n)}{\{(n+m-1)\beta/2+1\}^\beta_\kappa k!C^{{\beta}}_\kappa( I_m)}\\ \nonumber 
&=\sum_{k=0}^{\infty}\sum_{\kappa \in P^k_m} \frac{(n\beta/2)^\beta_\kappa 
\;C^{{\beta}}_\kappa(-\frac{\beta}{2}x\Sigma^{-1})}{\{(n+m-1)\beta/2+1\}^\beta_\kappa k!}
={{}_1F_1{}}^{(\beta;m)}\left(\frac{n\beta}{2};\frac{(n+m-1)\beta}{2}+1;-\frac{\beta}{2}x\Sigma^{-1}\right ).
\nonumber
\\
& = \mathrm{etr}\biggl(-\frac{\beta x}{2}\Sigma^{-1}\biggl){
 {}_1F_1{}}^{(\beta;m)}\left(\frac{(m-1)\beta}{2}+1;\frac{(n+m-1)\beta}{2}+1;\frac{\beta}{2}x\Sigma^{-1}\right ).
\label{betahyper}
\end{align}
The last equation \eqref{betahyper} follows by the $\beta$-Kummer relation ${{}_1F_1}^{(\beta;m)}(a;c;-X)=
\mathrm{etr}(-X){{}_1F_1}^{(\beta;m)}(c-a;c;X)$ given in D\'iaz-Garc\'ia~\cite{Garcia2014}. 
\end{proof}

The truncated distribution function of \eqref{eq:-cdf-l1-2} up to the $K$th degree is denoted by 
$F_K^{\beta}(x)$ in the simalar mannar of (\ref{eq:eq26}). 
Jack polynomials were computed by the improved algorithm of Hashiguchi et al.~\cite{Hashiguchi2000} and Hashiguchi and Niki~\cite{Hashiguchi2006}.
The original algorithm had been implemented for $\beta = 1$, but it is easy to improve for any $\beta$. 
As real application of wireless communication for $\beta=2$, 
we consider a multiple input and single output system (MISO).
Ratnarajah and Vaillancourt~\cite{Ratnarajah2005b} computed the channel capacity of a multiple inputs and 
multiple outputs (MIMO) under the population covariance matrix
 $\Sigma = \mathrm{diag}(1.81, 1.31, 0.69, 0.19)$, 
 where the numbers of target and transmit are $n=2$ and $m=4$, respectively. 
 Table~\ref{tbl:F2-100} shows the computation times under the above covariance and 
 MISO~($n=1$) for each dimension. 
 The last case of $m=4$ took 39h on a computer~(Mac Pro maxOS Catalina, ver. 10.15.5 with a 2.7 GHz Intel Xeon E5) and the used memory was 77GB.
 The channel capacity of MISO is provided by
\begin{align*}
C_1=\int_{0}^{\infty} \log_2\biggl(1+ \frac{\rho \ell_1}{m}\biggl) f(\ell_1) d\ell_1,
\end{align*}
where $\ell_1$ is the largest eigenvalue of a singular complex Wishart matrix and $\rho$ is the signal to noise ratio. 
Fig. \ref{fig:graph01} shows the capacity in nats vs.  signal to noise ratio for MISO under 
 $\Sigma=0.69I_3$ and $\mathrm{diag}(1.81, 1.31, 0.69)$.
 This figure is similar to Figure~2 of Ratnarajah and Vaillancourt~\cite{Ratnarajah2005a}. 
The two capacities are very close if the SNR is less than $11$ dB.

\begin{table}[H]
\begin{center}
\caption{Computation time of (\ref{eq:-cdf-l1-2}) for $n=1$, $\beta=2$, $K=100$ and each covariance matrix $\Sigma$ below, and the value of maximum probability that we computed.} 
\label{tbl:F2-100}
\begin{tabular}{clcc}
\multicolumn{1}{c}{$m$} & \multicolumn{1}{c}{$\Sigma$} & \multicolumn{1}{c}{Comp. time} 
& $\max_x F_K^\beta(x)$ \\ \hline
2 & $\mathrm{diag}(1.81,1.31)$ & 4.9s & 0.999 \\
3 & $\mathrm{diag}(1.81, 1.31, 0.69)$ & 648s & 0.999 \\
4 & $\mathrm{diag}(1.81, 1.31, 0.69, 0.19)$ & 39h & 0.973 \\ \hline
\end{tabular}
\end{center}
\end{table}
%%%%%%%%画像
\if0
\begin{figure}[H]
 \begin{center}
     \includegraphics[width=7cm]{./eps/ComFx(4,1).eps}
     \rlap{\raisebox{1.0ex}{\kern-0.0em{\small $x$}}}%
\rlap{\raisebox{28.0ex}{\kern-21.0em{\small $\Pr(\ell_1 < x) $}}}%
        \caption{$m=4$, $n=1$, $\Sigma=\mathrm{diag}(1.8090, 1.3090, 0.6910, 0.1910)$}\label{fig:graph01}
        \end{center}
  \label{fig:one}
\end{figure}
\fi
%%%%%%%%

\begin{figure} [H]
 \begin{center}
     \includegraphics[width=7cm]{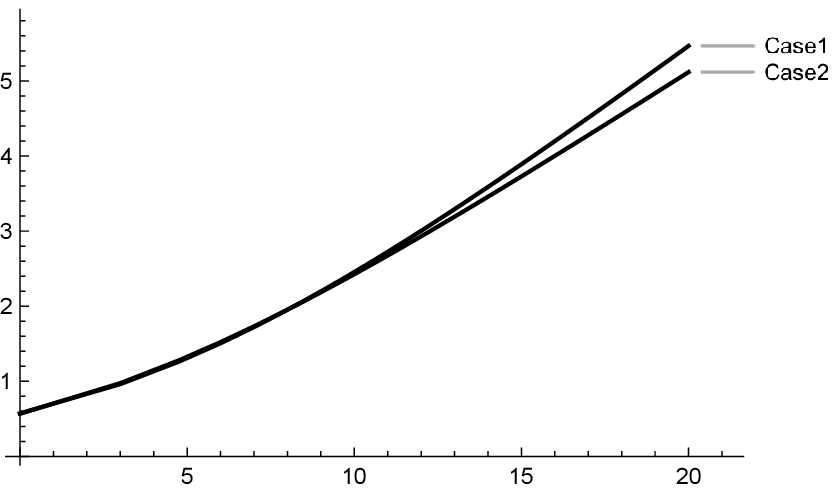}
     \rlap{\raisebox{26.0ex}{\kern-20.0em{\small $C_1$}}}%
     \rlap{\raisebox{1ex}{\kern-2.0em{\small \text{SNR}}}}%
        \caption{The channel capacity of MISO when uncorrelated and correlated cases.\\
 \text{Case1}:~ $\Sigma  = 0.69I_3$
\text{and }\text{Case2}:~$\Sigma =\mathrm{diag}(1.81, 1.31, 0.69).$
} \label{fig:graph01}
        \end{center}
  \label{fig:one}
\end{figure}

\section{Conclusion}
In this study, we considered the exact distributions of the largest eigenvalue of a singular beta-Wishart matrix. 
To derive the distributions of a singular beta-Wishart matrix, we defined heterogeneous hypergeometric functions with two matrix arguments. 
We provided a formula (\ref{splitting:beta}) that differs from Theorem~3 in D\'iaz-Garc\'ia~\cite{Garcia2013}. 
We applied the function (\ref{hyper0f0}) for $\beta=2$ to the channel capacity of MISO. 
Numerical computations were performed with a lower dimension $m$. 
Approximate distributions for distribution functions (\ref{hyper0f0}) and (\ref{betamax}) are still required. 
Furthermore, numerical computations for large sample sizes are planned as part of future work. 
Finally, the holonomic gradient method (HGM) proposed by Hashiguchi et al.~\cite{Hashiguchi2013,Hashiguchi2018} may also be applicable to such computations. 
%%%%%%%
\section*{Acknowledgements}
The authors would like to express their gratitude to the Editor-in-Chief, Associate Editors and the anonymous referees for careful reading and valuable suggestions. 
This work was supported by JSPS KAKENHI Grant Number 18K03428.
%%%%%%%
\section*{Appendix. Derivation of $f(\ell_1)$ for $m=2$}
 Let $W\sim \mathcal{W}_2(1,\Sigma)$; then the density function of $W$ is given as 
\begin{align*}
f(W)=\frac{1}{2\pi\sqrt{|\Sigma|}}\ell_1^{-1}\mathrm{etr}\biggl(-\frac{1}{2}\Sigma^{-1}W\biggl).
\end{align*}
Because $m=2$, we consider the case of $n=1$ and just one eigenvalue $\ell_1$.
The spectrum decomposition of $W$ is given as
\begin{align*}
W&=HLH^\top=\ell_1h_1h_1^\top=
\begin{pmatrix}
\ell_1  \cos^2{\theta}&\ell_1 \cos{\theta}\sin{\theta} \\
\ell_1 \sin{\theta}\cos{\theta}&\ell_1 \sin^2{\theta}
\end{pmatrix},
\end{align*}
where $H=(\boldsymbol{h}_1, \boldsymbol{h}_2)$, $\boldsymbol{h}_1=(\cos{\theta},\sin{\theta})^\top$, $\boldsymbol{h}_2=(-\sin{\theta},\cos{\theta})^\top$ and $L=\mathrm{diag}(\ell_1,0)$.
We note that $\boldsymbol{h}_2$ vanishes because of $\ell_2 = 0$.  
Thus, the density function $f(\ell_1,\theta)$ is given as
 \begin{align*}
f(\ell_1,\theta)
&=\frac{1}{2}\frac{1}{2\pi\sqrt{|\Sigma|}}\exp \biggl(-\frac{1}{2\lambda_1}\biggl)
 \exp \bigg\{-\frac{1}{2}\ell_1\biggl(\frac{1}{\lambda_2}-\frac{1}{\lambda_1}\sin^2 \theta \biggl)\biggl\},
\end{align*}
where $\lambda_1$ and $\lambda_2$ are the eigenvalues of $\Sigma$.
We integrate the density function of $f(\ell_1,\theta)$ as
\begin{align*}
\int_{0}^{2\pi}f(\ell_1,\theta)d\theta&=\frac{1}{4\pi \sqrt{|\Sigma|}}\exp \biggl(-\frac{1}{2\lambda_1}\biggl)\exp(a\sin^2\theta) d\theta=\frac{1}{4\pi \sqrt{|\Sigma|}}\exp \biggl(-\frac{1}{2\lambda_1}\biggl)\sum_{k=0}^{\infty} \int_{0}^{2\pi} \frac{a^k}{k!}(\sin^2\theta)^k d\theta ,
\end{align*}
where $a=-\frac{1}{2}\ell_1(\frac{1}{\lambda_2}-\frac{1}{\lambda_1})$. Hence we have
\begin{align*}
f(\ell_1)
=\frac{1}{4\pi \sqrt{|\Sigma|}}\exp \biggl(-\frac{1}{2\lambda_1}\ell_1\biggl) \sum_{k=0}^{\infty}2\pi \frac{a^k}{k!}\frac{(2k-1)!!}{2k!!}\nonumber = \frac{1}{2 \sqrt{|\Sigma|}}\exp \biggl(-\frac{1}{2\lambda_1}\ell_1\biggl){_1F_1}\biggl(\frac{1}{2};1;a\biggl). 
\end{align*}

   %%%%%%%%%%%%%%%%%%%%%%%%%%%%%%%%%%%%%%%%%%%%%%%%%%%%
   \bibliography{}

\end{document}